\documentclass[12pt]{amsart}
\usepackage{amsfonts}
\usepackage{graphicx}
\usepackage{amscd,hyperref}
\usepackage{hyperref}
\usepackage{enumitem}
\newtheorem{theorem}{Theorem}[section]
\theoremstyle{plain}
\newtheorem{acknowledgement}{Acknowledgement}

\newtheorem{corollary}{Corollary}[section]

\newtheorem{definition}{Definition}[section]
\newtheorem{example}{Example}

\newtheorem{lemma}{Lemma}[section]

\newtheorem{remark}{Remark}[section]

\numberwithin{equation}{section}

\begin{document}

\title[]{ An inequality concerning the growth bound of a discrete evolution family on a complex Banach space}
\author {Constantin Bu\c se, Donal O'Regan and Olivia Saierli}

\address{Politehnica University of Timisoara, Department of Mathematics, Piata Victoriei No. 2, 300006-Timisoara, Rom\^ ania }

\email{constantin.buse@upt.ro}

\email{buse1960@gmail.com}

\address{National University of Ireland, Galway, Ireland}
\email{donal.oregan@nuigalway.ie}

\address{Tibiscus University of Timisoara, Department of Computer Science and Applied Informatics, Str. Lasc\u ar Catargiu, No. 4-6,300559-Timisoara, Rom\^ania}
\email{saierli$\_$olivia@yahoo.com}
\maketitle

\begin{abstract} We prove that the uniform growth bound $\omega_0(\mathcal{U})$ of a discrete evolution
family $\mathcal{U}$ of bounded linear operators acting on a complex Banach space $X$ satisfies the
inequality
$$\omega_0(\mathcal{U})c_{\mathcal{U}}(\mathcal{X})\le -1;$$ here $c_{\mathcal{U}}(\mathcal{X})$ is the operator norm of a convolution operator which acts on a certain Banach space $\mathcal{X}$ of $X$-valued sequences.

\end{abstract}

\section{Notations, definitions and statement}

Let $X$ be a complex Banach space and let $\mathcal{L}(X)$ be the
Banach algebra of all bounded linear operators acting on $X.$ The
norm of $X$ and the operator norm on $\mathcal{L}(X)$ are denoted by
$\|\cdot\|.$ We use the classical notations $\mathbb{Z}_+$ and
$\mathbb{C}$  for the sets of nonnegative integers and of complex
scalars, respectively.
 As is well-known, the space $l^{\infty}(\mathbb{Z}_+, X)$ consisting by all bounded $X$-valued sequences becomes a Banach space when we endow it with the "sup" norm, i.e. $\|(f_n)\|_{\infty}=\sup\nolimits_{n\in\mathbb{Z}_+}\|f_n\|$. Let $l_0^\infty(\mathbb{Z}_+, X)$ be the subspace of $l^{\infty}(\mathbb{Z}_+, X)$ which consists of all sequences $(f_n)\in l^{\infty}(\mathbb{Z}_+, X)$ with $f_0=0.$ Also consider $c_{0}^{0}(\mathbb{Z}_+, X)$, the subspace of $l_{0}^{\infty}(\mathbb{Z}_+, X)$ consisting of all sequences $(f_n)$ having the property that $\lim\nolimits_{n\to\infty}f_n=0$. Obviously, $c_{0}^{0}(\mathbb{Z}_+, X)$ and $l_{0}^{\infty}(\mathbb{Z}_+, X)$ are closed subspaces of the Banach space $l^{\infty}(\mathbb{Z}_+, X).$ Now let  $1\le p<\infty.$ By $l_0^p(\mathbb{Z}_+, X)$ we denote the space of all $X$-valued sequences $f=(f_k)_{k\in\mathbb{Z}_+}$ having the property that $f_0=0$ and
 $$\|f\|_p:=\left(\sum\nolimits_{k=0}^\infty \|f_k\|^p\right)^{\frac{1}{p}}<\infty.$$
 Obviously, $(l_0^p(\mathbb{Z}_+, X), \|\cdot\|_p)$ is a Banach space.

 \begin{definition}\label{def} A family $\mathcal{U}:=\{U(n, m): (n, m)\in\mathbb{Z}_+\times\mathbb{Z}_+, n\ge m\}\subset\mathcal{L}(X)$ is
  called a discrete evolution family if it satisfies the properties: $U(m,m) = I$
  and $U(m, n) = U(m, p)U(p, n)$ for all nonnegative integers $m\geq p\geq n.$\end{definition}

   Here $I$ denotes the identity operator on $X$.

Let us denote by $\Delta$ the set of all pairs of nonnegative integers $(n, m),$ so that $n\ge m$ and
let $\Omega(\mathcal{U})$ be the set of all real numbers $\omega$ such that
\begin{equation}\label{(1)}\sup\limits_{(n, m)\in\Delta}e^{-\omega(n-m)}\|U(n, m)\|:=M_{\omega}<\infty.\end{equation}

Throughout the paper we assume that $\Omega(\mathcal{U})$ is a non-empty set, i.e. the
family $\mathcal{U}$ is exponentially bounded. The uniform growth
bound of $\mathcal{U}$, denoted by $\omega_0(\mathcal{U}),$ is the
infimum of $\Omega(\mathcal{U})$.

A typical example which provides a discrete evolution family is
presented next.

\begin{example} Let $\mathcal{A}:=\{A_n: n\in\mathbb{Z}_+\}$ be a family of bounded linear operators acting on a Banach space $X.$ The discrete evolution family associated to the family $\mathcal{A}$ is the two parameters family $\mathcal{U}_{\mathcal{A}}:=\{U_{\mathcal{A}}(m,n): m\geq n\in\mathbb{Z}_+\}\subset\mathcal{L}(X)$ given by
 $$U_{\mathcal{A}}(m,n):=\left\{\begin{array}{lc}
                    A_{m-1}A_{m-2}\cdots A_n, \quad m>n \\
                    I,\quad m=n.
                  \end{array}
 \right.$$
Obviously, the family $\mathcal{U}_{\mathcal{A}}$ is an evolution family in the sense of Definition \ref{def}.  Moreover, every evolution family $\mathcal{U}:=\{U(n, m): (n, m)\in\mathbb{Z}_+\times\mathbb{Z}_+, n\ge m\}$ comes in this way. It is enough to set $A_n=U(n+1, n)$ in order to see this.

 For any $X$-valued sequence $f=(f_n)_{n\in \mathbb{Z}_+}$ we consider the discrete inhomogeneous Cauchy Problem
\begin{equation}\label{discrete-CP}
\left\{\begin{array}{lc}
         x_{n+1}=A_nx_n+f_{n+1}, \quad n=0, 1,\ldots \\
         x_0=0.
       \end{array}
\right.
\end{equation}
Obviously, the solution of (\ref{discrete-CP}) is the sequences $(x_n)$, given by
$$x_n=\sum\limits_{k=0}^n U_{\mathcal{A}}(n, k)f_k,\quad n\in\mathbb{Z}_+.$$

\end{example}

Let $\mathcal{M}:=\{c_{0}^{0}(\mathbb{Z}_+, X),
l^{\infty}_0(\mathbb{Z}_+, X), l_{0}^{p}(\mathbb{Z}_+, X)\},$
$\mathcal{X}\in\mathcal{M}$ and let $\mathcal{U}:=\{U(n, m): (n, m)\in\mathbb{Z}_+\times\mathbb{Z}_+, n\ge m\}$ be an exponentially bounded discrete evolution family.

For each $j\in\mathbb{Z}_+$ and each sequence
$f=(f_n)\in\mathcal{X}$ consider the linear operator
$\mathcal{T}_{\mathcal{X}}(j)$ given by
$$(\mathcal{T}_{\mathcal{X}}(j)f)(k):=\left\{\begin{array}{clc}
                                 U(k, k-j)f_{k-j}, & \mbox{ for all } (k, j)\in\Delta \\
                                 0, & \mbox{ otherwise.}
                               \end{array}
\right.$$ Since the family $\mathcal{U}$ is exponentially bounded, $\mathcal{T}_{\mathcal{X}}(j)$ is well defined and acts on $\mathcal{X}$. The
family ${\bf
\mathcal{T}_{\mathcal{X}}}:=\{\mathcal{T}_{\mathcal{X}}(j)\}_{j\in\mathbb{Z}_+}$
is a discrete semigroup, i.e. $\mathcal{T}_{\mathcal{X}}(0)$ is the
identity operator on $\mathcal{X}$ and
$\mathcal{T}_{\mathcal{X}}(j+k)=\mathcal{T}_{\mathcal{X}}(j)\circ\mathcal{T}_{\mathcal{X}}(k)$
for all nonnegative integers $j$ and $k.$ It is called the evolution
semigroup associated to the discrete family $\mathcal{U}$ on
$\mathcal{X}$.

Let $T$ in $\mathcal{L}(X) $ be a single operator. In the following,
$\rho(T)$ denotes the resolvent set of $T$, i.e.  the set of all
complex scalars $z$ for which $zI-T$ has a bounded inverse in
$\mathcal{L}(X)$. Also $\sigma(T):=\mathbb{C}\setminus\rho(T)$
denotes the spectrum of the operator $T.$ As is well-known the
spectrum of $T$ is a compact and non-empty set. The spectral radius
of $T$, denoted by $r(T)$, is defined as $r(T):=\sup\{|z|:
z\in\sigma(T)\}$. It is well known (Gelfand spectral radius theorem,
$1941$) that
$$r(T)=\lim\limits_{n\rightarrow\infty}\|T^n\|^{\frac{1}{n}}.$$

Obviously, this yields
\begin{equation}\label{(2)} \ln(r(T))=\lim\limits_{n\to\infty}\frac{\ln\|T^n\|}{n}=\omega_0(\{U(n, m):=T^{n-m}: (n, m)\in\Delta\}).\end{equation}

Since $\|T^n\|\leq\|T\|^n,$  $r(T)\le\|T\|$ and hence $\sigma(T)$ is
a subset of $\{z\in\mathbb{C}:|z|\leq\|T\|\}$. For each
$z\in\rho(T), R(z, T):=\left(zI-T\right)^{-1}$ denotes the resolvent
operator of $T$. It is well-known that for every $z\in\rho(T),$ one
has
\begin{equation}\label{50}\|R(z, T)\|\cdot \mbox{ dist }(z, \sigma(T))\ge 1.\end{equation} In
particular, if $z_n\in\rho(T)$ and $z_n\to z\in\sigma(T)$ then $\|R(z_n, T)\|\to\infty$ as $n\to\infty.$

The series
$\left(\sum\nolimits_{n\geq 0}\frac{T^n}{z^{n+1}}\right)$
is absolutely convergent on $\{|z|>r(T)\}$ and its sum is given by
\begin{equation}\label{3}
    \sum\nolimits_{n=0}^\infty\frac{T^n}{z^{n+1}}= \frac{1}{z}\sum\nolimits_{n=0}^\infty\left(\frac{T}{z}\right)^n=\frac{1}{z}\cdot\left(I-\frac{T}{z}\right)^{-1}= R(z, T),
  \end{equation}
for every $z\in\mathbb{C}$ with $|z|>\|T\|.$

The next Lemma (whose proof is in Section 2) connects the uniform
growth bound of an exponentially bounded evolution family
$\mathcal{U}$ and the spectral radius of
$\mathcal{T}_{\mathcal{X}}(1).$
\begin{lemma}\label{lemma1} Let $\mathcal{U}$ be an exponentially bounded evolution family as given
above, $\mathcal{X}\in\mathcal{M}$ and let $\mathcal{T}_{\mathcal{X}}$ be the discrete evolution semigroup associated to $\mathcal{U}$ on $\mathcal{X}.$ Then

  \begin{equation}\label{(10)}\omega_0(\mathcal{U})=\ln r(\mathcal{T}_{\mathcal{X}}(1)).\end{equation}

\end{lemma}

The "convolution" operator $\mathcal{K}_{\mathcal{X}}:D(\mathcal{K}_{\mathcal{X}})\subset \mathcal{X}\to \mathcal{X},$ associated to the discrete family $\mathcal{U},$
is defined by
\begin{equation} D(\mathcal{K}_{\mathcal{X}}):=\{f\in\mathcal{X}: U\ast f\in\mathcal{X}\}\end{equation} where
\begin{equation}   (U\ast f)(k):=\sum\nolimits_{j=0}^k U(k, j)f_j,\quad f=(f_n)\in \mathcal{X}, k\in\mathbb{Z}_+.
\end{equation}

\begin{theorem}\label{thm-1.1.} Let $\mathcal{U}, \mathcal{X}$, $\mathcal\mathcal{T}_{\mathcal{X}}$ and $\mathcal{K}_{\mathcal{X}}$ be as above. The following statements are equivalent.
\begin{itemize}
  \item[{\bf 1.}] The family $\mathcal{U}$ is uniformly exponentially stable, i.e. its uniform growth bound is negative.
  \item[{\bf 2.}] The evolution semigroup $\mathcal{T}_{\mathcal{X}}$ associated to the family $\mathcal{U}$ on $\mathcal{X}$ is uniformly
exponentially stable, i.e. $r(\mathcal\mathcal{T}_{\mathcal{X}}(1))$ is less than one.
  \item[{\bf 3.}] For each $f\in \mathcal{X},$  $U\ast f$ belongs to $\mathcal{X}$.
  \item[{\bf 4.}] The linear operator $f\mapsto\mathcal{K}_{\mathcal{X}}(f)$ is bounded on $\mathcal{X}.$
\end{itemize}

\end{theorem}
For further details, counterparts or different versions of the above
result we refer the reader to \cite{[CD]},\cite{[BCDS]},\cite{[BKRT]},\cite{[KRZ]} and the references therein.

In the continuous case, results like the previous one are well-known. For further details we refer the reader to \cite{[CL]}, \cite{[CLMR]}, \cite{[D]}, \cite{[MRS]} and the references therein.

The proof of Theorem \ref{thm-1.1.} for $\mathcal{X}\in\{c_0^0(\mathbb{Z}_+, X), l_0^{\infty}(\mathbb{Z}_+, X)\}$ is the same as in \cite[Thm.~3.4]{[BKRT]}. We mention that the 4th statement in Theorem \ref{thm-1.1.} is not contained in the statement of \cite[Thm.~3.4]{[BKRT]}, but its equivalence with the first three statements is established in the proof.

However for the discrete $l_0^p(\mathbb{Z}_+, X)$-version of the above theorem we could not find a  reference in the literature. The proof of the present version is similar to that given in \cite[Thm.~3.4]{[BKRT]} so we present part of the argument when $\mathcal{X}=l^p_0(\mathbb{Z}_+, X)$.

  To prove that the first statement implies the third one, let $N$ and $\nu$ be two positive constants such that $\|U(n, m)\|\le Ne^{-\nu(n-m)}$ for every $(n, m)\in\Delta$ and let $f\in l_0^p(\mathbb{Z}_+, X).$ Thus
$$ \begin{array}{crc}
     \|U*f\|_p^p &\le& N^p\sum\nolimits_{n=0}^\infty e^{-\nu np}\sum\nolimits_{k=0}^n e^{\nu kp}\|f_k\|^p  \\
    & \le & N^p\sum\nolimits_{k=0}^\infty e^{\nu kp}\|f_k\|^p \sum\nolimits_{n=k}^\infty e^{-\nu n p} \\
    & \le & \frac{N^pe^{\nu p}}{e^{\nu p}-1}\|f\|_p^p.
   \end{array}$$
For the proof of ${\bf 3.}\Rightarrow {\bf 4.}$ we can argue as in the proof of the first step in \cite[Thm.~3.3]{[BKRT]} and we mention  the well-known fact that convergence in $l^p_0(\mathbb{Z}_+, X)$ implies convergence on coordinates.

 Now, we prove that the last statement implies the first one. Since $\mathcal{K}_{\mathcal{X}}$ is bounded, there exists a positive constant $c_p$ such that

  \begin{equation}\label{61}\|\mathcal{K}_{\mathcal{X}}f\|_p\le c_p\|f\|_p \mbox{ for all } f\in l_0^p(\mathbb{Z}_+, X).\end{equation}
  Let $j\ge 1$ be an integer and $x\in X.$ Let $f=(f_n)\in l_0^p(\mathbb{Z}_+, X)$ with $f_j=x$ and $f_k=0$ whenever $k$ is different of $j,$ and set $U(n, j)=0$ when $n<j.$ Thus inequality (\ref{61}) yields
  $$\sum\nolimits_{n=j}^{\infty}\|U(n, j)x\|^p\le c_p^p\|x\|^p$$ and from the discrete version of the well-known Datko theorem it follows that the family $\mathcal{U}$ is uniformly exponentially stable; see for example \cite[Thm.~3.4]{[PPP]}.

  The connection with the second statement can be made by following the proof of \cite[Thm.~3.4]{[BKRT]}. Finally, we mention that Lemma 3.1 and Theorem 3.5 from \cite{[BKRT]} remain valid, with the same proof if $l_0^p(\mathbb{Z}_+, X)$ or $l_0^{\infty}(\mathbb{Z}_+, X)$ replaces $c_{00}(\mathbb{Z}_+, X).$ Also, we mention that the assumption $x(0)=0$ in Lemma 3.1 from \cite{[BKRT]} is essential. This is the reason why we consider spaces of sequences having first entry equal to $0,$

 When the family $\mathcal{U}$ is uniformly exponentially stable, we
 let
 $$c_\mathcal{U}(\mathcal{X}):=\|\mathcal{K}_{\mathcal{X}}\|_{\mathcal{L}(\mathcal{X})}=\sup\limits_{\|f\|_{\mathcal{X}}\leq 1}\|U\ast f\|_{\mathcal{X}}.$$

\begin{theorem}\label{thm-1.2} Let $\mathcal{X}\in\mathcal{M}$ and let $\mathcal{U}$ be a uniformly exponentially stable evolution family acting on $X.$ Then the following three statements hold true.
\begin{itemize}
\item{\bf (i)} The following inequality occurs:
\begin{equation}\label{thm-1.2-eq}
\omega_0(\mathcal{U})\cdot c_\mathcal{U}(\mathcal{X})\le -1.\end{equation}
\item{\bf (ii)} The resolvent set $\rho(\mathcal{T}_{\mathcal{X}}(1))$ contains the set $$\pi:=\left\{|z|> 1-\frac{1}{c_{\mathcal{U}}(\mathcal{X})}\right\}.$$
\item{\bf (iii)} The resolvent operator satisfies the estimate
\begin{equation}\sup\limits_{|z|\ge 1}\|R(z, \mathcal{T}_{\mathcal{X}}(1))\|\le  c_{\mathcal{U}}(\mathcal{X}).\end{equation}
\end{itemize}
\end{theorem}

If the discrete evolution family $\mathcal{U}$ satisfies the
convolution condition
$$U(n, j)=U(n-j, 0)\mbox{ for all } (n, j)\in\Delta,$$
then with $T=U(1, 0)$ we have that $T^n=U(n, 0).$ In this case, the
convolution operator is defined by
$$(\mathcal{S}_{\mathcal{X}}f)(n):=(T\ast f)(n)=\sum\nolimits_{j=0}^n T^{n-j}f_j,\quad f=(f_j)\in\mathcal{X}.$$
Obviously, $\mathcal{S}_{\mathcal{X}}$ acts on $\mathcal{X}$ and it is a bounded linear operator on $\mathcal{X}$ provided that $r(T)<1.$ For further details concerning similar results for strongly continuous semigroups see for example \cite{[vanNe]}. In this particular case, for each pair $(n, m)\in\Delta$, we have that $U(n, m)=T^{n-m}$. Moreover, $c_{\mathcal{U}}(\mathcal{X})=\|\mathcal{S}_{\mathcal{X}}\|$ and $r(\mathcal{T}_{\mathcal{X}}(1))=r(T).$  The above Theorem \ref{thm-1.2} reads as
\begin{corollary}\label{cor-1} Let $T$ be a single operator in $\mathcal{L}(X)$ such that $r(T)<1$ and let $\mathcal{X}\in\mathcal{M}.$ Then
\begin{equation}\label{30}-1\geq \|\mathcal{S}_{\mathcal{X}}\|_{\mathcal{L}(\mathcal{X})}\ln(r(T)).\end{equation}

\end{corollary}

A natural question to ask is if the inequality (\ref{30}) is sharp.
The next example shows that it can be arbitrarily tight.
\begin{example} Let $X=\mathbb{C}$, $T:=\gamma\in(0, 1)$ and $\mathcal{X}=c_0^0(\mathbb{Z}_+, X).$ Thus $\|T^n\|=\gamma^n$ and $r(T)=\|T\|=\gamma<1.$
Also
$$\|\mathcal{S}_{\mathcal{X}}\|=\sup\nolimits_{n\in\mathbb{Z}_+}\sum\nolimits_{k=0}^n\gamma^k=\frac{1}{1-\gamma}$$
and the inequality (\ref{30}) becomes
$$\ln(\gamma)\cdot\frac{1}{1-\gamma}\le-1.$$ The equality is attained for $\gamma\to
1$ (l'H\^{o}pital's rule).

\end{example}

  We note that the above Theorem \ref{thm-1.1.} does not provide a negative number $\sigma$ such
  that $\omega_0(\mathcal{U})$ is less than $\sigma$ while our result does this.

  As is well known, if $T\in\mathcal{L}(X)$ and
  $$\sum\nolimits_{n=0}^\infty\|T^nx\|<\infty,\quad \forall x\in X$$
  then $r(T)<1$; see \cite{[Gluck]} for updated results of this type. For  comprehensive information on this subject we refer the reader to \cite{[vN1]}. In some sense, this result can be improved to
  \begin{corollary}\label{cor-2}
  Let $T\in\mathcal{L}(X)$ such that  $\sum_{n=0}^\infty\|T^n\|:=u_1(T)<\infty$. Then for each $1\le p<\infty$, one has
\begin{equation}\label{60}
\ln(r(T))\cdot u_1(T)\le \|\mathcal{S}_{l_0^p(\mathbb{Z}_+,X)}\|_{\mathcal{L}(l_0^p(\mathbb{Z}_+,X))}\ln(r(T))\leq -1.
\end{equation}
  \end{corollary}
%---------------------------------------------------------------
%
%
%
%
%---------------------------------------------------------------
\section{Proofs}

We start this section with the proof of Lemma \ref{lemma1}. We already stated that the discrete evolution family $\mathcal{U}$ is uniformly exponentially stable if and only if $r(\mathcal{T}_{\mathcal{X}}(1)$ is less than $1.$ Lemma \ref{lemma1} can be derived from this by a simple scaling argument, as was already done; see for example \cite[Thm.~3.5]{[BKRT]}. However, for completeness, we  present here a direct proof  using only the definitions stated above. Let
$\omega\in\Omega(\mathcal{U})$, $\mathcal{X}\in\mathcal{M}$ and $f\in\mathcal{X}$. After an obvious calculation, we get
$$\|\mathcal{T}_{\mathcal{X}}(m)f\|_\mathcal{X}\leq M_\omega e^{\omega m}\|f\|_\mathcal{X},\quad \forall m\in\mathbb{Z}_+,$$
where $M_\omega$ is defined in (\ref{(1)}). Therefore,
$$\frac{\ln\left(\|\mathcal{T}_{\mathcal{X}}(1)^m\|\right)}{m}\leq \frac{\ln M_\omega}{m}+\omega,\quad \forall m\geq 1.$$
Based on (\ref{(2)}),  the previous inequality yields
$\ln(r(\mathcal{T}_{\mathcal{X}}(1)))\leq\omega$, which produces
\begin{equation}\label{(3)}
\ln(r(\mathcal{T}_{\mathcal{X}}(1)))\leq\omega_0(\mathcal{U}).
\end{equation}

In order to establish the reverse inequality in (\ref{(3)}), let
$j\in\mathbb{Z}_+$, $j\geq 1$, $x\in X$, $x\neq 0$ and set
$$f_k=\left\{\begin{array}{clc}
               x,& \mbox{ if } k=j\\
               0,& \mbox{ otherwise}.
             \end{array}
\right.$$
Obviously, $f=(f_k)$ belongs to $\mathcal{X}.$ Let $\nu>\omega_0(\mathcal{T}_\mathcal{X})=\ln(r(\mathcal{T}_\mathcal{X}(1)))$ (see (\ref{(2)})). Thus there exist $K_{\nu} \ge 1$ such that
$$\|\mathcal{T}_{\mathcal{X}}(j)\|_{\mathcal{L}(\mathcal{X})}\le K_{\nu}e^{\nu j}, \quad\forall j\in\mathbb{Z}_+.$$ Therefore, for every $n\in\mathbb{Z}_+,$ one has
\begin{equation}\label{(4)}
\|U(n+j, j)x\|=\|(\mathcal{T}_{\mathcal{X}}(n)f)(n+j)\|\leq \|\mathcal{T}_{\mathcal{X}}(n)\|_{\mathcal{L}(\mathcal{X})}\|x\|\le K_{\nu}e^{\nu }\|x\|.
\end{equation}
On the other hand
\begin{equation}\label{(5)}
\|U(n+1, 0)x\|=\leq  \|\mathcal{T}_\mathcal{X}(n)\|_{\mathcal{L}(\mathcal{X})}\|U(1, 0)\|\|x\|\le K_{\nu}e^{\nu }\|U(1, 0)x\|,
\end{equation}
where the inequality (\ref{(4)}) was used. From (\ref{(4)}) and
(\ref{(5)}) we have that $\omega_0(\mathcal{U})\leq\ln(r(\mathcal{T}_\mathcal{X}(1))),$ which completes the proof.

\emph{Proof of Theorem \ref{thm-1.2}.}

 For every $z\in \mathbb{C}$, $|z|=1$, $n\in\mathbb{Z}_+$ and $f\in\mathcal{X}$, one has
$$\begin{array}{cclc}
    \left[R(z,\mathcal{T}_\mathcal{X}(1))f\right](n) & = &\sum\limits_{k=0}^\infty\frac{(\mathcal{T}_\mathcal{X}(k)f)(n)}{z^{k+1}}\\
                                                     &= & \sum\limits_{k=0}^n\frac{U(n, n-k)f_{n-k}}{z^{k+1}}\\
                                                     &= & \frac{1}{z^{n+1}}\sum\limits_{j=0}^nU(n, j)(z^jf_{j})\\
                                                     &= & \frac{1}{z^{n+1}}(U\ast g)(n),\\
  \end{array}$$
where $g_j:=z^jf_j$ for all $j\in\mathbb{Z}_+$ and $g=(g_j)$. Clearly, $f\in\mathcal{X}$ if and only if $g\in\mathcal{X}$ and, in addition $\|g\|_\mathcal{X}=\|f\|_\mathcal{X}$. Hence
\begin{equation}\|R(z,\mathcal{T}_\mathcal{X}(1))f\|_\mathcal{X}=\|U\ast g\|_\mathcal{X}\leq \|
\mathcal{K}_{\mathcal{X}}\|_{\mathcal{L}(\mathcal{X})}\|g\|_\mathcal{X}= c_\mathcal{U}(\mathcal{X})\|f\|_\mathcal{X},\end{equation} i.e.
\begin{equation}\label{Rc}
\|R(z,\mathcal{T}_\mathcal{X}(1))\|\leq c_\mathcal{U}(\mathcal{X}), \quad \forall z\in\mathbb{C}, |z|=1.
\end{equation}

Thus for any $z\in\mathbb{C}$ with $|z|=1$ we have that
\begin{equation}\label{Tc}  c_\mathcal{U}(\mathcal{X})\ge \| R(z,\mathcal{T}_\mathcal{X}(1))\|\ge\frac{1}{1-r(\mathcal{T}_\mathcal{X}(1))},  \end{equation}
where different counterparts of the result in \cite[Thm.~3.5]{[BKRT]} and (\ref{50}) with $\mathcal{T}_\mathcal{X}(1))$ instead of $T$, was used.

Now, (\ref{thm-1.2-eq}) is a consequence of the elementary inequality
\begin{equation}\frac{1}{1-r}\ge \frac{-1}{\ln(r)}, \quad r\in(0, 1). \end{equation} Thus statement {\bf (i)} is settled.

On the other hand, (\ref{Tc}) can be written in the form
\begin{equation}r(\mathcal{T}_{\mathcal{X}}(1))\le 1-\frac{1}{c_{\mathcal{U}}(\mathcal{X}}\end{equation}
which readily yields {\bf (ii)}, as well.

Since the Neuman series expansion of the resolvent shows that $$\|R(z, \mathcal{T}_{\mathcal{X}}(1))\|\to 0 \mbox{ as }|z|\to\infty,$$ assertion {\bf (iii)} follows from (\ref{Rc}) and from the Phragm$\acute{e}$n-Lindel\"of theorem.

\begin{remark} We used \cite[Thm.~3.5]{[BKRT]} and its counterparts in the proof of Theorem \ref{thm-1.2}, although it is not needed at all to derive (\ref{Tc}). Indeed, let us choose $\lambda\in\sigma(\mathcal{T}_{\mathcal{X}}(1)$ with $|\lambda|=r(\mathcal{T}_{\mathcal{X}}(1))$ and a complex number $z$ with $|z|=1$ and $\arg(z)=\arg(\lambda).$ Then one readily obtains

\begin{equation}  c_\mathcal{U}(\mathcal{X})\ge \| R(z,\mathcal{T}_\mathcal{X}(1))\|\ge\frac{1}{\mbox{dist}(z, \sigma(\mathcal{T}_{\mathcal{X})}(1)}=\frac{1}{|\lambda-z|}=\frac{1}{1-r(\mathcal{T}_\mathcal{X}(1))}.\end{equation} We thank the referee who brought our attention concerning this important fact.

\end{remark}

\vspace{0.5cm}
\emph{Proof of Corollary \ref{cor-2}.} We divide the proof into two parts by considering the cases $p=1$ and $p>1$ separately. Let $f=(f_k)\in\l_0^1(\mathbb{Z}_+, X)$. Then
\begin{equation}\begin{array}{cclc}
    \|T\ast f\|_1 & = & \sum\nolimits_{n=0}^\infty\left\|\sum\nolimits_{k=0}^n T^{n-k}f_k\right\| \\
                  & \leq &  \sum\nolimits_{n=0}^\infty\sum\nolimits_{k=0}^\infty 1_{\{0, 1, \ldots, n\}}(k)\|T^{n-k}\|\|f_k\| \\
                  & = & \sum\nolimits_{k=0}^\infty\|f_k\|\sum\nolimits_{n=k}^\infty \|T^{n-k}\| \\
                  &= & u_1(T)\cdot \|f\|_1,
  \end{array}
\end{equation}
which yields
$$\|\mathcal{S}_{l_0^1(\mathbb{Z}_+,X)}\|_{\mathcal{L}(l_0^1(\mathbb{Z}_+, X))}\leq  u_1(T).$$

As is usual, $1_B$ denotes the characteristic function of the set $B\subset \mathbb{Z}_+.$

Multiplying the above inequality by $\ln(r(T))$ and using  Corollary \ref{cor-1} (with $l_0^1(\mathbb{Z}_+,X)$ instead of $\mathcal{X}$) we get
$$-1\geq \|\mathcal{S}_{l_0^1(\mathbb{Z}_+,X)}\|_{\mathcal{L}(l_0^1(\mathbb{Z}_+,X))}\ln(r(T))\geq u_1(T)\ln(r(T)).$$

 When $p\in(1, \infty)$, let $f=(f_k)\in l^p_0(\mathbb{Z}_+, X)$ and $h=(h_k)\in l^q_0(\mathbb{Z}_+, X^\ast)$, with $\frac{1}{p}+\frac{1}{q}=1$. Then
\begin{equation}
\begin{array}{clc}
  \left\|\sum_{k=0}^\infty h_k(T\ast f)(k)\right\| & \leq\sum_{k=0}^\infty  \left\|h_k\right\| \left\|(T\ast f)(k)\right\| \\
              & \leq  \sum_{k=0}^\infty  \left\|h_k\right\| \sum_{j=0}^k \left\|T^j\right\|\left\|f_{k-j}\right\| \\
              & = \sum_{k=0}^\infty \sum_{j=0}^\infty  \left\|h_k\right\| 1_{\{0,\cdots, k\}}(j)\left\|T^j\right\|\left\|f_{k-j}\right\|\\
              & = \sum_{j=0}^\infty \left\|T^j\right\| \sum_{k=j}^\infty  \left\|h_k\right\| \left\|f_{k-j}\right\|. \end{array}
\end{equation}

 Using H\"older's inequality, we get
 \begin{equation}
\begin{array}{clc}
  \left\|\sum_{k=0}^\infty h_k(T\ast f)(k)\right\|
               &\leq \sum_{j=0}^\infty\left\|T^j\right\| \left(\sum_{k=j}^\infty  \left\|h_k\right\|^q\right)^{1/q}\|f\|_p \\
              &\leq \|h\|_q\|f\|_p\sum_{j=0}^\infty\|T^j\|\\
              &= \|h\|_q\|f\|_p u_1(T).
\end{array}
\end{equation}
Thus \begin{equation}\label{eq21}\|T\ast
f\|_p=\sup\limits_{\|h\|_q\leq 1}\left\|\sum\nolimits_{k=0}^\infty h_k(T\ast
f)(k)\right\|\leq \|f\|_p u_1(T),\end{equation} and
hence $u_1( T)\ge \|f\mapsto T\ast
f\|_{\mathcal{L}(l^p_0(\mathbb{Z}_+, X))}.$ The assertion follows by using Corollary \ref{cor-1} and
taking into account that $\ln(r(T))$ is negative.

\begin{acknowledgement}

The authors would like to thank the referees for their help and suggestions in improving this paper.
 
\end{acknowledgement}

%--------------------------------------------------------------
%
%
%
%
%---------------------------------------------------------------

\end{document}